\documentclass[11pt]{article}
\usepackage[utf8]{inputenc}
\usepackage[
papersize={5.5in, 8.5in},
left=0.35in,
right=0.35in,
top=0.35in,
bottom=0.5in]{geometry}
\usepackage{graphicx}
\usepackage{amsthm}
\usepackage{rotating}
\usepackage{amsfonts,amssymb,amsmath}
\usepackage[usenames]{color}
\usepackage[T1]{fontenc}

\usepackage{array}
\usepackage{float}
\usepackage{booktabs}
\usepackage[dvipsnames,usenames,table,xcdraw]{xcolor}

\graphicspath{ {images/} }

\title{\bf 
On lower bounds of the order \\of $k$-chromatic unit distance graphs 
}

\author{\normalsize {Aubrey D.N.J. de Grey; \textcolor[rgb]{0.2,0.6,0.4}{Jaan Parts (jaan\_parts@mail.ru)} }}
\date{\normalsize {Mountain View, California, USA; \textcolor[rgb]{0.2,0.6,0.4}{Kazan, Tatarstan, Russia}}}

\begin{document}

\maketitle

\pagestyle{empty}
\thispagestyle{empty}

\begin{abstract}
Here we give refined numerical values for the minimum number of vertices of $k$-chromatic unit distance graphs in the Euclidean plane.
\end{abstract}

\section{Background}

In the previous issue of \textit{Geombinatorics}, an article \cite{gwyn} was published, in which Haydn Gwyn and Jacob Stavrianos obtained new estimates for the \textit{minimum} number of vertices $v_k$ and edges $e_k$ for arbitrary $k$-chromatic\footnote{
For clarity, we emphasize that here $k$ means the number of colors in the proper vertex coloring of the \textit{graph} and differs by one from similar designations in \cite{gwyn} and other publications mentioned below, where it is associated with the \textit{tiling} of the plane.}
\textit{unit distance graphs} in the plane for $k=5$ and $k=6$. 

It is known that $v_3=3$ and $v_4=7$ (provided by the unit triangle and the Moser graph). For $k\ge 5$, exact values of $v_k$ are not known. Initial \textit{lower bounds} $v_5>12$, $v_7>6197$ were obtained by Dan Pritikin \cite{pri}. For $k=5$, the finite \textit{upper bounds} are\footnote{In \cite{par1} the corresponding 509-vertex graph has 2442 edges, but is not edge-critical, which allows us to reduce $e_5$. We were able to discard 36 edges, but we didn't perform an exhaustive search, so further improvements are possible.}
$v_5\le 509$, $e_5\le 2406$, see \cite{par1}. 

Pritikin's approach is based on finding a proper tiling of the plane using $k-1$ colors, which covers most of the plane, so that the average \textit{density} $\delta$ of the \textit{uncolored} area is minimal. Thus, $k-1$ colors are enough for any unit distance graph with at most $\lfloor 1/\delta \rfloor$ vertices, and a $k$-chromatic graph must have at least one more vertex.

In \cite{par2} we slightly improved the lower bounds to $v_6>24$ and $v_7>6992$ by modifying the tiling construction.

The main idea of Gwyn-Stavrianos' approach is to first estimate the minimum number of edges $e_k$, and then proceed to the number of vertices $v_k$ using known relations. To do this, a tiling of the plane in $k-1$ colors is found, minimizing the \textit{probability} $p_{k-1}$ of the formation of a \textit{mono-chromatic edge} (such that both vertices have the same color). The value $1/p_{k-1}$ gives a lower bound on the number of edges $e_k$.

For $k=5$, \cite{gwyn} uses a tiling by regular hexagons of four colors with a diameter of about 1.1335. For $k=6$, it is overlaid with a set of disks of unit diameter, located at a unit distance from each other with an average density $\delta=\pi/(8\sqrt3)$ and colored in the fifth color, and a lemma is used that allows passing from $p_4$ to $p_5$:\; $p_k\le(1-2\delta) \;p_{k-1}$.

According to the calculations of Gwyn and Stavrianos, $e_5\ge 98$, $v_5\ge 22$, $e_6\ge 180$, $v_6\ge 32$.
In this note, we offer several refinements that allow us to slightly increase these numerical values.

\section{Improvements}

First, note that in \cite{gwyn}, to go from the number of edges $e_k$ to the number of vertices $v_k$, the relation $e_k<v_k^{3/2}$ was used, which was obtained by Paul Erdős in \cite{erd}. Péter Ágoston and Dömötör Pálvölgyi recently obtained tighter bounds \cite{ago}:
$e_k\le \sqrt[3]{29/4}\;v_k^{4/3}$. In the range $20\le v_k< 521$, better estimates can be derived from the formula\footnote{
For completeness, note that according to \cite{ago} for $1\le v_k\le15$ the minimum number of edges $e_k$ is known exactly: 0, 1, 3, 5, 7, 9, 12, 14, 18, 20, 23, 27, 30, 33, and 37; while for $16\le v_k\le19$ the upper bounds on $e_k$ are 42, 47, 52, and 57. 
For the next $20\le v_k< 60$, using the formula, we get the following upper bounds on $e_k$: 63, 68, 72, 77, 82, 87, 92, 97, 102, 108, 113, 119, 124, 130, 136, 142, 148, 154, 160, 166, 172, 179, 185, 192, 198, 205, 212, 218, 225, 232, 239,
246, 254, 261, 268, 276, 283, 291, 298, 306.}
\cite{ago}: $v_k \ge 2\,r-11+24/v_k\,+$ $(1-s)\binom{\lfloor r \rfloor}{2} + s \binom{\lceil r \rceil}{2}$, where $r=2e_k/v_k$, $s=r-\lfloor r \rfloor$.
Applying these relations to $e_k$, we get a noticeable improvement in the $v_k$ estimates.

In \cite{gwyn} in the case of $k=6$, the plane was partially covered by disks with density $\delta=\pi/(8\sqrt3)\approx 0.226725$. Instead, one can use 
Croft's tiling \cite{cro} with rounded 12-gons and a density of about $\delta\approx 0.229365$, which slightly improves the $e_6$ estimate.
Also note that in \cite{gwyn} rounding down was used to get the $e_k$ estimates. It is more correct to use rounding up $e_k\ge\lceil 1/p_{k-1} \rceil$, which gives an increase by one more.

In \cite{gwyn}, a statistical method was used to find the value of $p_k$: a repeating section $S$ of the tiling was selected, on which a random unit edge was repeatedly superimposed, defined by the coordinates $\sigma=(x, y)$ of its first vertex and by the orientation angle $\phi$ of its second one, and the average value of the binary function $M_k(\sigma,\phi)$ was determined, which takes the value 1 or 0 depending on whether an edge is mono-chromatic or not. 
Instead, we took an algebraic method, calculating the integral 
$p_k=\frac{1}{2\pi S}\int_S d\sigma\int_{0}^{2\pi}d\phi \,M_k(\sigma,\phi)$
in the \texttt{Mathematica} package, which allowed us to refine the optimal tiling parameters. However, we encountered some problems with the numerical integration function \texttt{NIntegrate} in the general case, and were forced to limit ourselves to relatively simple tilings.

As a result, we get the bounds: $e_5\ge 99$, $v_5\ge 28$, $e_6\ge 182$, $v_6\ge 42$.

We also tried Gwyn-Stavrianos' approach for the case $k=7$, and obtained the following estimates for a tiling close to Pritikin's construction using pentagonal tiles: $e_7\ge 232646$, $v_7\ge 6456$. (For the construction using hexagonal tiles, we were able to get only $e_7\ge 69451$ and $v_7\ge 2608$.) This is better than \cite{pri}, but falls short of \cite{par2}. However, there is still some hope of beating the latter one by using tiles with curved borders (the so-called "wavy edges" used in \cite{par2}).

\end{document}